\documentclass[12pt,oneside]{amsart}
\usepackage{amssymb,amsmath}
\makeatletter
\def\@cite#1#2{{\m@th\upshape\bfseries%
[{#1\if@tempswa{\m@th\upshape\mdseries, #2}\fi}]}}
\makeatother
\theoremstyle{plain}
\newtheorem{thm}{Theorem}[section]
\newtheorem{lem}[thm]{Lemma}
\newtheorem{cor}[thm]{Corollary}

\theoremstyle{definition}
\newtheorem{rem}[thm]{Remark}

\newtheorem{defn}[thm]{Definition}

\newtheorem{egs}[thm]{Examples}

\newtheorem{prob}[thm]{Problem}
\newcommand{\Prf}{\noindent\textbf{Proof.\ }}
\newcommand{\bx}{\strut\hfill$\blacksquare$\medbreak}

\newcommand{\ca}{\mathrm{C}^*}

\DeclareMathOperator*{\wotalg}{\textsc{wot}--{\rm Alg}}

\newcommand{\wot}{\textsc{wot}}

\newcommand{\bbC}{{\mathbb{C}}}

\newcommand{\bbF}{{\mathbb{F}}}

\newcommand{\bbQ}{{\mathbb{Q}}}

\newcommand{\bbT}{{\mathbb{T}}}
\newcommand{\bbZ}{{\mathbb{Z}}}
 \newcommand{\A}{{\mathcal{A}}}
 \newcommand{\B}{{\mathcal{B}}}

\renewcommand{\H}{{\mathcal{H}}}
 \newcommand{\I}{{\mathcal{I}}}
 \newcommand{\J}{{\mathcal{J}}}

 \newcommand{\M}{{\mathcal{M}}}
 
\renewcommand{\O}{{\mathcal{O}}}

\newcommand{\upchi}{{\raise.35ex\hbox{$\chi$}}}
\newcommand{\fA}{{\mathfrak{A}}}

\newcommand{\fL}{{\mathfrak{L}}}

\newcommand{\fR}{{\mathfrak{R}}}

\newcommand{\qand}{\quad\text{and}\quad}

\newcommand{\qwhere}{\quad\text{where}\quad}

\newcommand{\Alg}{\operatorname{Alg}}

\newcommand{\dist}{\operatorname{dist}}

\newcommand{\Lat}{\operatorname{Lat}}

\newcommand{\sat}{\operatorname{sat}}
\newcommand{\fngee}{\bbF^+\!(G)}

\newcommand{\flgee}{\fL_G}
\newcommand{\frgee}{\fR_G}

\title{\normalsize\bf \rule{0pt}{3.5cm} PARTLY FREE ALGEBRAS  \\
FROM DIRECTED GRAPHS}

\author{DAVID W. KRIBS AND STEPHEN C. POWER}

\thanks{first author partially supported by a Canadian NSERC
Post-doctoral Fellowship.}  
\date{} 
\begin{document}
\maketitle
\thispagestyle{empty}

\baselineskip=12pt
\begin{quote}
We say that a nonselfadjoint operator algebra is {\it partly free} if it 
contains a free semigroup algebra. Motivation for such 
algebras occurs in the setting of what we call {\it free semigroupoid 
algebras}. These are the weak operator topology closed algebras generated 
by the left regular representations of semigroupoids associated with 
finite or countable directed graphs. We expand our analysis of partly
free algebras from previous work and obtain a graph-theoretic characterization of when a free
semigroupoid algebra
with countable  graph is partly free. 
This analysis carries over to norm closed {\it quiver algebras}. 
We also discuss new examples for the countable graph case. 
\end{quote} 

\vskip 1truecm
\baselineskip=15pt

%

Every finite or countable directed graph $G$ recursively generates a Fock space Hilbert space
and a family of partial isometries. 
These operators are of Cuntz-Krieger-Toeplitz type and  also arise through the left regular 
representations of  free semigroupoids determined by directed graphs. 
This was initially discovered by Muhly 
\cite{M1}, 
and, in the case of finite graphs, more recent work with Solel \cite{MS1} 
considered the norm
closed 
algebras generated by these representations, which they called {\it quiver 
algebras}. In \cite{KP1}, we developed a structure theory for the weak  
operator topology closed algebras $\flgee$ generated by the left regular  
representations coming from both finite and countable directed graphs; we 
called these algebras {\it free semigroupoid algebras}. In doing so, we 
found a 
unifying framework for a number of classes of algebras which appear in the 
literature, including; noncommutative analytic Toeplitz algebras $\fL_n$  
\cite{AP,DP1,DP2,Kfactor,Pop_fact,Pop_beur} (the prototypical free 
semigroup 
algebras), the classical analytic Toeplitz 
algebra $H^\infty$ \cite{Douglas_text,Hoffman_text}, and certain finite dimensional 
digraph algebras \cite{KP1}.  
But this approach gives rise to a diverse collection of new examples 
which include finite dimensional algebras, algebras with free behaviour, algebras which can be 
represented as matrix function algebras, and examples which mix these possibilities.  The
general theme of our work in 
\cite{KP1} was a marriage of simple 
graph-theoretic properties with properties of the operator algebra. 
Furthermore, our technical analyses were chiefly  
spatial in nature; for instance, we 
proved the graph is a complete unitary invariant of both the free semigroupoid algebra and the
quiver algebra.

In the next section we give a short introduction to free
semigroupoid algebras and discuss a number of examples.  The second
section contains an expanded analysis of 
a subclass called {\it partly free algebras}, which are characterized by containment of a copy of a
free 
semigroup algebra. Specifically, we extend  analysis from \cite{KP1} to 
the case of countable graphs, obtaining a new graph-theoretic description 
of when a free semigroupoid algebra is partly free. 
In fact, we find that it is unusual for an algebra coming from a countable graph to not 
be partly free. 
We present a number of new illustrative examples for the countable graph 
case. 
Our analysis works equally well for quiver algebras, and we obtain a 
graph-theoretic condition for when these norm closed algebras are partly 
free.

\section{Free Semigroupoid Algebras}\label{S:structure}

Let $G$ be a finite or countable directed graph, with edge set $E(G)$ and 
vertex set $V(G)$. Let $\fngee$ be the {\it free semigroupoid} determined 
by 
$G$; that is, $\fngee$ consists of the vertices, which act as units, and 
all allowable finite paths in $G$, with the natural operations of 
concatenation of allowable paths. Given a path $w$ in $\fngee$ we write 
$w=ywx$ when the initial and final vertices of $w$ are, respectively, $x$ 
and $y$. Let $\H_G = \ell^2(\fngee)$ be the Hilbert space with 
orthonormal basis indexed by elements of $\fngee$. For each edge $e\in 
E(G)$ and vertex $x\in V(G)$, we may define partial isometries and 
projections on $\H_G$ by 
the following actions on basis vectors:
\[
L_e\xi_w = \left\{ \begin{array}{cl}
\xi_{ew} & \mbox{if $ew\in\fngee$} \\
0 & \mbox{otherwise}
\end{array}\right.
\]
and
\[
L_x\xi_w = \left\{ \begin{array}{cl}
\xi_{xw}=\xi_w & \mbox{if $w = xw\in\fngee$} \\
0 & \mbox{otherwise}
\end{array}\right.
\]

These operators may thus be regarded as `partial creation operators' 
acting on a generalized
Fock space Hilbert space. There is an equivalent tree perspective  \cite{KP1}
which gives an appealing visual interpretation of the actions of 
these operators.
The family $\{L_e,L_x\}$ also arises through the left regular representation
$\lambda_G: \fngee \rightarrow \B(\H_G)$, with $\lambda_G(e) = L_e$, and 
$\lambda_G(x) =
L_x$. 
The associated {\it free semigroupoid algebra} is the weak operator 
topology closed algebra generated by this family, 
\begin{eqnarray*}
\flgee &=& \wotalg\,\, \{ L_e,L_x:e\in E(G), x\in V(G) \} \\
&=& \wotalg \,\,\{ \lambda_G(w) : w\in\fngee \}.
\end{eqnarray*}

\begin{rem}
In the case of finite graphs, Muhly and Solel \cite{M1,MS1} considered the 
norm closed algebras $\A_G$ generated by such a family, calling them {\it 
quiver algebras}. For both finite and countable graphs, we considered the classification problem
for $\A_G$ in \cite{KP1},  and in Section~2.1 we derive  partly free conditions  for $\A_G$. 
Recently the $\ca$-algebras generated by
families of partial isometries associated with directed graphs have been studied heavily. The set
of generators for these algebras are sometimes referred to as 
Cuntz-Krieger $E$-families (for instance see  
\cite{BHRS,Ephrem,Kumjian1,Kumjian2}). 
On the other hand, the generators of free semigroupoid algebras are of Cuntz-Krieger-Toeplitz
type in the sense that the   
$\ca$-algebra generated by a family $\{ L_e \}$ is generally the extension of a         
Cuntz-Krieger    $\ca$-algebra by the  compact operators. 
\end{rem}


There is also a right regular representation 
$\rho_G:\fngee\rightarrow
\B(\H_G)$ determined by $G$, which yields partial isometries $\rho_G(w) \equiv
R_{w^\prime}$  for
$w\in\fngee$ acting on $\H_G$ by the equations $R_{w^\prime}\xi_v = 
\xi_{vw}$, where
$w^\prime$ is the word $w$ in reverse order, with similar conditions. The
corresponding algebra is 
\begin{eqnarray*}
\frgee &=& \wotalg \,\, \{R_e,R_x:e\in E(G), x\in V(G) \} \\
&=& \wotalg \,\, \{ \rho_G(w) : w\in\fngee \}.
\end{eqnarray*}
Given edges $e,f\in E(G)$, observe that $L_eR_f \xi_w = \xi_{ewf} = R_fL_e\xi_w$, for all
$w\in\fngee$, so that $L_e R_f = R_f L_e$, and similarly for the vertex projections. In fact,
the commutant of $\frgee$ coincides with $\flgee = \frgee^\prime$ \cite{KP1}. Further,  the
commutant of 
$\flgee$ coincides with $\frgee$, and thus $\flgee$ is its own second commutant. 

A useful ingredient in the proof is the observation that the algebras $\frgee$ and $\fL_{G^t}$ are
naturally unitarily equivalent, where $G^t$ is the {\it transpose directed 
graph} obtained from $G$ by reversing
the directions of all edges. An important technical device obtained here is 
the existence of
Fourier expansions for elements of $\flgee$. Specifically, if $A\in\flgee$ and $x\in V(G)$, then
$A\xi_x = \sum_{w=wx} a_w \xi_w$ with $a_w\in\bbC$, and 
Cesaro type sums associated with the formal sum $\sum_{w\in\fngee} a_w L_w$, converge in
the strong operator topology to $A$. We write, $A\sim \sum_{w\in\fngee} a_w L_w$. As a
notational convenience, 
for projections determined by vertices $x\in V(G)$, we put $P_x \equiv L_x$ 
and $Q_x \equiv R_x$.

In \cite{KP1} we proved that $G$ is a complete unitary invariant of both $\flgee$ and $\A_G$.
Thus, different directed graphs really do yield different algebras. We finish this section by
discussing a number of examples from simple graphs. 

\begin{egs}
${\mathbf (i)}$
The algebra generated by the graph with a single vertex $x$ and single loop edge $e = xex$ is
unitarily
equivalent to the classical analytic Toeplitz algebra $H^\infty$
\cite{Douglas_text,Hoffman_text}. Indeed, the Hilbert space in this case may be naturally
identified with the Hardy space $H^2$ of the unit disc, and under this 
identification $L_e$ is easily seen to be
unitarily equivalent to the unilateral shift $U_+$. 

${\mathbf (ii)}$ 
The noncommutative analytic Toeplitz algebras $\fL_n$, $n\geq 2$
\cite{AP,DKP,DP1,DP2,Kfactor,Pop_fact,Pop_beur}, the fundamental examples 
of  free semigroup algebras,
arise from the graphs with a single vertex and $n$ distinct loop edges. 
For instance, in the case $n =2$ with loop edges $e=xex \neq f=xfx$, the Hilbert space is
identified with unrestricted 2-variable Fock space $\H_2$. The operators $L_e, L_f$ are
equivalent to the natural creation operators on $\H_2$, also known as the 
Cuntz-Toeplitz
isometries. Further, $P_x = I$, and thus $\flgee \simeq \fL_2$. 

${\mathbf (iii)}$
As an example of a simple matrix function algebra, we may consider the graph $G$ with
vertices $x, y$ and edges $e=xex$, $f=yfx$. Then $\flgee$ is generated by 
$\{ L_e, L_f , P_x ,
P_y \}$. If we make the natural identifications $\H_G = P_x \H_G \oplus P_y \H_G \simeq H^2
\oplus H^2$ (respecting word length), then 
\[
L_e \simeq \left[ 
\begin{matrix} 
U_+ & 0 \\
0 & 0 
\end{matrix} \right]
\quad
L_f \simeq \left[ 
\begin{matrix} 
0 & 0 \\
U_+ & 0 
\end{matrix} \right]
\quad
P_x \simeq \left[ 
\begin{matrix} 
I & 0 \\
0 & 0 
\end{matrix} \right]
\quad
P_y \simeq \left[ 
\begin{matrix} 
0 & 0 \\
0 & I 
\end{matrix} \right].
\]
Thus, $\flgee$ is seen to be unitarily equivalent to 
\[
\flgee \simeq \left[
\begin{matrix} 
H^\infty & 0 \\
H^\infty_0 & \bbC I  
\end{matrix}\right] 
\]
where $H^\infty_0$ is the subalgebra of $H^\infty$ functions $h$ with $h(0) = 0$.

${\mathbf (iv)}$
By simply adding a directed edge $g = xgy$ to the previous graph, we obtain a very different
algebra $\fL_{G^\prime}$. In fact, $\fL_{G^\prime}$ is a unitally partly free algebra in the
sense of the next section because it contains isometries with mutually orthogonal ranges; for
instance, $U = L_e^2 + L_f L_g$ and $V = L_e L_g + L_g L_e$ are isometries which satisfy
$U^* V =0$. 
 
${\mathbf (v)}$
If $G$ is a finite graph with no directed cycles, then the Fock space $\H_G$ is                    
finite-dimensional and so too is $\flgee$. As an example, consider the graph with three vertices
and two edges, labelled $x_1$, $x_2$, $x_3$, $e$, $f$
where $e = x_2ex_1$, $f= x_3fx_1$. Then the Fock space is spanned by the vectors $\{
\xi_{x_1}, \xi_{x_2}, \xi_{x_3}, \xi_e, \xi_f \}$ and with this basis
the general operator
$
X = 
\alpha L_{x_1} + \beta L_{x_2} + \gamma L_{x_3} + \lambda L_e + \mu L_f
$
in $\flgee$ is represented by the matrix
\[
X \simeq 
\begin{bmatrix}
\alpha & & & & \\
 & \beta & & & \\
 & & \gamma & & \\
\lambda & & & \beta & \\
\mu & & & & \gamma
\end{bmatrix}.
\]
Algebraically, $\flgee$ is isometrically isomorphic to the so-called digraph algebra in
$\M_3(\bbC)$ consisting of the matrices
\[
\begin{bmatrix} 
\alpha & 0 & 0 \\ \lambda & \beta & 0 \\ \mu & 0 & \gamma 
\end{bmatrix}.
\]
Recall that a digraph algebra $\A(H)$ is a unital subalgebra of $\M_n(\bbC)$ which is spanned
by some of the standard matrix units of $\M_n(\bbC)$. The graph $H$ is transitive and reflexive
and is such that the edges of $H$ naturally label the relevant matrix units.

${\mathbf (vi)}$
Let $n\geq 1$ and consider the {\it cycle graph} $C_n$ which has $n$
vertices $x_1, \ldots, x_n$  and $n$ edges $e_n =x_1 e_n x_n$ and $e_k =x_{k+1} e_k x_k$ for
$k=0,\ldots, n-1$. 
The {\it cycle algebra} $\fL_{C_n}$ may be identified with the $\wot$-closed semicrossed
product
$\bbC^n \times_\beta^\sigma \bbZ_+$ associated with the cyclic shift automorphism $\beta$ of
$\bbC^n$ \cite{DeAPe}.  To see this identify $L_{x_i}\H_G$ with $H^2$ for each $i$ in the
natural way
(respecting word length). Then $\H_G = L_{x_1}\H_G \oplus \ldots \oplus L_{x_n}\H_G \simeq
\bbC^n \otimes H^2$ and the operator $\alpha_1 L_{e_1} + \ldots + \alpha_n L_{e_n}$ is
identified with the operator matrix 
\[
\left[
\begin{matrix}
0 & & & & \alpha_n T_z \\
 \alpha_1 T_z & 0 & & &  \\
 & \alpha_2 T_z & 0 & &  \\
 & & \ddots & \ddots &  \\
 & & & \alpha_{n-1}T_z & 0  
\end{matrix}\right].
\]

Writing $H^\infty (z^n)$ for the subalgebra of $H^\infty$ arising from functions of the form
$h(z^n)$ with $h$ in $H^\infty$, the algebra $\fL_{C_n}$ is readily identified with the matrix
function algebra 
\[
\left[
\begin{matrix}
H^\infty (z^n) & z^{n-1}H^\infty (z^n)& \hdots & z H^\infty (z^n) \\
 z H^\infty (z^n) & H^\infty (z^n) & & \vdots  \\
\vdots &   & \ddots &  \\
z^{n-1} H^\infty (z^n) & \hdots &  & H^\infty (z^n)  
\end{matrix}\right].
\]
This in turn is identifiable with the crossed product above. In fact, this matrix function algebra is
the $\wot$-closed variant of the matrix function algebra $\B_n$ of De Alba and Peters
\cite{DeAPe} for the norm closed semicrossed product $\bbC^n \times_\beta^\sigma \bbZ_+$.
Such
identifications are the Toeplitz versions of the identification of the graph $\ca$-algebra of $C_n$
with $\M_n ({\mathrm C} (\bbT))$. 

\end{egs}

\section{Partly Free Algebras}\label{S:partfree}

We say that a $\wot$-closed operator algebra
$\fA$ is {\it partly free} if it contains the free semigroup algebra $\fL_2$ as a subalgebra in the
sense of the following definition.  
\begin{defn}\label{partlyfreedefn}
A $\wot$-closed algebra $\fA$ is {\it partly free} if there is an inclusion map $\fL_2
\hookrightarrow \fA$ which is the restriction of an injection between the generated von
Neumann algebras. If the map can be chosen to be unital, then $\fA$ is said to be {\it unitally
partly free}.
\end{defn}
These notions parallel the requirement that a $\ca$-algebra contain
the Cuntz algebra $\O_2$, or that a discrete group contain a free group.
Theorems~\ref{weak} and \ref{strong} determine when the
algebras $\flgee$ are partly free and unitally partly free for both finite and countable directed
graphs. 
We require the following structural result on partial isometries in $\flgee$. 

\begin{lem}\label{standardform}
\cite{KP1}
The initial projections of partial isometries $V$ in $\flgee$ are sums of projections from $\{ P_x
:x\in V(G) \}$. Specifically, 
\[
V^* V = \sum_{x\in \I} P_x  
\qwhere \I = \big\{x\in V(G) : V\xi_x \neq 0\big\}.
\]
\end{lem}

\noindent{In} fact, we proved much more in \cite{KP1}. All partial 
isometries $V$ in $\flgee$ satisfy a
standard form, written as $V = \sum_{x\in \I} \oplus L_{\eta_x}$, where $\{\eta_x\}_{x\in\I}$
are unit {\it wandering} vectors for $\frgee$ supported on distinct $Q_x\H_G$, and hence the
initial projections
satisfy $L_{\eta_x}^* L_{\eta_x} = P_x$. 


The cycle algebras $\fL_{C_n}$, $1 \leq n < \infty$, arise in the general theory, as does the
algebra $\fL_{C_\infty}$ generated by an infinite graph analogue of $C_n$. Let $C_\infty$ be
the directed graph  
with vertex set $\{ x_k : k \geq 1\}$ indexed by the natural numbers, and edge set 
$\{e_k = x_{k+1} e_k x_k :k\geq 1\}$. In the finite graph case, the cycle algebras $\fL_{C_n}$ 
are the 
key examples of 
free semigroupoid algebras which are not partly free. However, 
$\fL_{C_\infty}$ is unitally partly free, and this gives an indication of  
how pervasive these algebras are in the countable graph case. 

\begin{lem}\label{infinitecyclegraph}
The algebras $\fL_{C_n}$, $1 \leq n < \infty$, are not partly free. 
However, $\fL_{C_\infty}$ is unitally partly free. 
\end{lem}

\Prf
We proved this result in \cite{KP1} for $1\leq n < \infty$ by showing 
$\fL_{C_n}$ does not contain 
pairs of partial
isometries $U$, $V$ which satisfy condition $(iii)$ of
Theorem~\ref{weak}. For $n = \infty$, we shall construct a pair of 
isometries $U$, $V$ in 
$\fL_{C_\infty}$ with orthogonal ranges. 
For $k\geq 1$, let $u_k$ and $v_k$ be the unique finite paths in $\bbF^+(C_\infty)$ with $u_k =
x_{2k}
u_k x_k$ and $v_k = x_{2k+ 1} v_k x_k$.
Then we may define $U$ and $V$ by
\[
U = \sum_{k \geq 1} \oplus L_{u_k}   
\qand 
V = \sum_{k \geq 1} \oplus L_{v_k},   
\]
where the sums converge in the strong operator topology. 
A unital injection of $\fL_2$ into $\fL_{C_\infty}$ is then defined   by 
mapping the two
generators of $\fL_2$ to $U$ and $V$.  
\bx

In fact it turns out  there are two ways in which an algebra $\flgee$ 
can be partly free; namely $G$ must either contain a double-cycle or a 
proper infinite (directed) path, in the sense of the following definition. 

\begin{defn}
A {\it cycle} at vertex $x$ is a path of edges $w=xwx$ for which only the 
initial edge has source vertex $x$. 
The graph $G$ contains a {\it double-cycle} if there are distinct cycles 
$w_i = xw_ix$,
$i=1,2$, at some vertex $x$ in $G$. 
By a {\it proper infinite (directed) path} in $G$, we mean an infinite 
word $\omega = 
\cdots e_{i_3}
e_{i_2}e_{i_1}$ in the edges of $G$ such that no edges are repeated 
and every finite segment 
 corresponds to an allowable finite directed path in $G$. 
\end{defn}

A little thought shows that the only way a graph $G$ fails to have either 
a double-cycle or such a proper infinite path is if either there are no 
infinite paths at all, or, if every infinite path in the graph becomes 
periodic; in this case say that $G$ has the periodic path property. If, 
on the other hand, $G$ does not have the periodic path property, we shall 
say that $G$ has the {\it aperiodic path property} meaning that there 
exists an aperiodic infinite path. This in turn means that there exists a 
proper infinite path or a double-cycle. 

We define the {\it saturation} at a vertex $x$ in
$G$ to be the set, $\sat (x)$, that consists of $x$, together with all 
finite 
and infinite paths which start at $x$,
and all vertices that are final vertices for paths starting at $x$.  
Consider a graph which is a downward directed tree in which every branch 
has finite 
length and different branches may intersect. Also perform  
the following surgery: add arbitrary cycles to some vertices without 
creating double-cycles. Note there can be infinite branchings. This 
`aperiodic looped tree 
graph' does not have the aperiodic path property and furthermore such a 
graph is a `fairly typical' graph without that property. In fact,  
$G$ does not have  the aperiodic path property if and only if the 
saturation graph 
for each vertex has this kind of structure. 

For the sake of brevity in the next proof we shall assume the finite 
vertex case \cite{KP1}. 




\begin{thm}\label{weak}
The following assertions are equivalent for a finite or countable directed graph $G$:
\begin{itemize}
\item[$(i)$] $G$ has the aperiodic path property.  
\item[$(ii)$] $\flgee$ is partly free.
\item[$(iii)$] There are nonzero partial isometries $U$, $V$ in $\flgee$ with
\[
U^*U = V^*V, \,\,\,\, UU^*\leq U^*U, \,\,\,\, VV^*\leq V^*V, \,\,\,\, U^* V =0.
\]
\end{itemize}  
\end{thm}

\Prf
Since von Neumann algebra isomorphisms are spatial, condition $(iii)$ is a reformulation of
$(ii)$. Thus it suffices 
to prove the equivalence of $(i)$ and $(iii)$.

To see $(iii)\Rightarrow (i)$, first suppose $U$ and $V$ are partial 
isometries in $\flgee$ satisfying $(iii)$ with $U^*U = P_x = V^*V$ for 
some $x\in V(G)$. 
Then there is at least one cycle
$w=xwx$ over $x$, for otherwise $U=P_xUP_x=0$ or $U=P_x$ and $0 = P_x U 
=V$. Suppose $w$ is the only cycle
over $x$. Then the compression algebra $P_x \flgee P_x$ forms
the $\wot$-closed subalgebra of $\flgee$ consisting of elements $X\in\flgee$ with Fourier
expansions of the form $X\sim \sum_{k\geq 0} a_{k} L_{w^k}$. Let $\H_x = P_x \H_G$.
Then $P_x \flgee|_{\H_x}$ is 
evidently unitarily equivalent to $\fL_{C_1}\simeq H^\infty$. But $P_x \flgee|_{\H_x}$
contains a
pair of non-zero partial isometries $U|_{\H_x}=P_x U|_{\H_x}$, $V|_{\H_x}=P_x V|_{\H_x}$
satisfying 
condition $(iii)$, and this
contradicts Lemma~\ref{infinitecyclegraph}. Thus, we deduce the existence of at least two
cycles over $x$, and hence $G$ contains a double-cycle. 

For the general case, let $U$ and $V$ be partial isometries in $\flgee$ 
satisfying
$(iii)$ with $U^*U=V^*V=\sum_{x\in\I} P_x$ as in Lemma~\ref{standardform}. 
Suppose, by way of 
contradiction, that $(i)$ fails so that $G$ has the structure discussed 
above. The vertices $x$ in the index set $\I$ are scattered over this 
graph. Plainly, either there exists such an $x$ with no other $y$ in $\I$ 
in the saturation of $x$ or there is a cycle $C$ containing several such 
edges such that the saturation of $C$ contains no further vertices $x$ in 
$\I$. In the first case we can consider the compression of the given $U$ 
and $V$ to $\H_x$, and in the second case the compression to $\H_C = 
\sum_{x\in V(C)} P_x\H_G$. 
Observe that the saturations of $x$ and $C$ may be infinite. However, the 
conditions in $(iii)$ allow us to see that these compressions are indeed 
partial isometries because the hypothesis implies that the compression is 
the same (modulo a zero summand) as the restriction to the subspace 
$\flgee \H_x$ or $\flgee \H_C$. 
Now we get our contradiction by arguing as in the previous paragraph. 
Indeed, the first case is trivial and in the second case the compression 
of the algebra is unitarily equivalent to a cycle graph algebra 
$\fL_{C_n}$.

Towards the implication $(i)\Rightarrow (iii)$, suppose $\omega$ is a 
proper infinite
path in $G$. Let $\J$ be the set of all vertices in $G$ for which the saturation at
each $x\in\J$ contains part, and hence an entire tail, of $\omega$. Let $\{x_k,y_k\}\rightarrow
\{z_k\}$ be a two-to-one
map from the countable set $\J$ onto itself, such that the vertices $x_k,y_k\in\J$ belong to the
saturation of $z_k\in\J$. Thus, there are distinct finite paths $u_k = x_ku_kz_k$, $v_k =
y_kv_kz_k$. By design, $U=\sum_k \oplus L_{u_k}$, $V=\sum_k \oplus L_{v_k}$ are partial
isometries in $\flgee$ with orthogonal ranges, initial projection $U^*U = P_\J = V^*V$, and
$U = P_\J UP_\J$, $V = P_\J VP_\J$. 

Next, suppose distinct cycles $w_1$, $w_2$ form a double-cycle over a vertex $x$ in
$G$. Let $\J$ be the set of all vertices for which the saturation at each vertex in $\J$ contains
this double-cycle. Let us enumerate the (possibly finite) vertices of $\J$ as $\{x_k\}_k$, 
and let $r_k = xr_kx_k$ be a path from $x_k $ to $x$. For each vertex $x_k$ in $\J$ choose,
without
repeating any choices, two paths $u_k$, $v_k$ amongst the set $\{ w_1^m w_2:m\geq 1\}$.
Then
again by design, $U=\sum_k \oplus L_{u_kr_k}$, $V=\sum_k \oplus L_{v_kr_k}$ are partial
isometries in $\flgee$ with orthogonal ranges, initial projection $U^*U = P_\J = V^*V$, and
$U = P_\J UP_\J$, $V = P_\J VP_\J$. Therefore, in both cases we have shown that
$(i)\Rightarrow (iii)$, and this completes the proof. 
\bx

It follows from this result that the free semigroupoid algebras coming from  countable graphs are
typically partly free.
Indeed, the condition on a countable graph $G$ 
which forces $\flgee$ to not be partly free    is quite
restrictive as outlined in the previous discussion.  Below we discuss a 
number of examples. 

The following is the unital version of the previous theorem.   
We shall say $G$ has the
{\it uniform aperiodic path property} if the saturation at every vertex 
includes an aperiodic infinite path.

\begin{thm}\label{strong}
The following assertions are equivalent for a finite or countable directed graph $G$:
\begin{itemize}
\item[$(i)$] $G$ has the uniform aperiodic path property. 
\item[$(ii)$] $\flgee$ is unitally partly free.
\item[$(iii)$] There are  isometries $U$, $V$ in $\flgee$ with
\[
U^* V =0.
\]
\end{itemize}  
\end{thm}

\Prf
Once again, condition $(iii)$ is a restatement of $(ii)$, so it suffices 
to prove the equivalence of $(i)$ and $(iii)$. 

For $(iii)\Rightarrow (i)$, the proof of Theorem~\ref{weak} can be 
adapted to show that the
saturation at every vertex in the index set of vertices $\I$, determining 
the initial projection for
the partial isometries $U$, $V$, includes a double-cycle or a 
proper infinite path. Hence we may apply this argument in the current 
case with $\I =V(G)$ as we are
dealing with isometries here, $U^*U=V^*V=I=\sum_{x\in V(G)} P_x$, and it follows that
$G$ has the uniform aperiodic path property. 

To see $(i)\Rightarrow (iii)$, consider the last two paragraphs in the proof of
Theorem~\ref{weak}. As $G$ satisfies the uniform aperiodic path property, 
it follows that we may
decompose the vertex set for $G$ into disjoint subsets $V(G) = \cup_i \J_i$, where each $\J_i$
is obtained as in one of these two cases; double-cycles, or proper  infinite 
paths. In either case we can define partial isometries $U_i$,
$V_i$ in $\flgee$ with orthogonal ranges, $U_i  = P_{\J_i} U_i P_{\J_i}$, $V_i
 = P_{\J_i} V_i P_{\J_i}$, and $U^*_i U_i = P_{\J_i} = V^*_iV_i$. Thus,  the operators $U =
\sum_i \oplus U_i$, $V = \sum_i \oplus V_i$ are isometries in $\flgee$
with mutually orthogonal ranges.
\bx

We next add to the short list of known hyper-reflexive algebras 
\cite{Arv_nest,Berc,Chr,Dav_Toep,DP1} by extending our result from 
\cite{KP1} to the case of countable graphs. Given an operator algebra
$\fA$, a measure of the
distance to $\fA$ is given by
\[
\beta_{\fA} (X) = \sup_{L\in\Lat\fA} ||P_L^\perp X P_L||,
\]
where $P_L$ is the projection onto the subspace $L$ and $\Lat \fA$ is the lattice of invariant
subspaces for $\fA$. Evidently,
$\beta_\fA (X) \leq \dist (X,\fA)$, and the algebra $\fA$ is said to be
{\it hyper-reflexive} if there is a constant $C$ such that $\dist(X,\fA)
\leq C \beta_\fA (X) $ for all $X$. 

The  free semigroup algebras $\fL_n$ were proved to be hyper-reflexive by 
Davidson ($n=1$ \cite{Dav_Toep}) and 
Davidson and Pitts ($n\geq 2$ \cite{DP1}). Furthermore,  motivated by the $\fL_n$ case 
Bercovici \cite{Berc} proved   an algebra is hyper-reflexive
with distant constant no greater than $3$ whenever its commutant contains
a pair of isometries with orthogonal ranges. 

\begin{cor}\label{hyperreflex}
Let $G$ be a finite or countable directed graph such that the transpose 
graph $G^t$
satisfies the uniform aperiodic path property; equivalently, $\fL_{G^t}$ is 
unitally partly free.  
Then $\flgee$ is hyper-reflexive with distant constant at most 3.
\end{cor}

\Prf
This is a direct consequence  of Bercovici's result  \cite{Berc}
since $\flgee^\prime = \frgee \simeq \fL_{G^t}$. 
\bx

\begin{prob}
Is $\flgee$ hyper-reflexive for every directed graph $G$?
\end{prob}


We next present some simple examples for the countable graph 
case. The focus will be on the new aspect discovered here; the relevance 
of proper infinite paths.  

\begin{egs}
${\mathbf (i)}$ Let $G$ be the directed graph with vertices $\{ x_k : k\in\bbZ\}$ indexed by the
integers 
and directed edges $\{e_k = x_{k+1} e_k x_k : k\in\bbZ\}$. We could also add (possibly
infinite) directed paths $w_k = x_k
w_k$. Then every vertex saturation in $G$ contains a proper infinite tail 
$\omega_k =
\cdots e_{k+1} e_k$ for some $k\in\bbZ$. Thus $G$ satisfies the 
uniform aperiodic path property,
and  $\flgee$ is unitally partly free.
 
Consider the interesting special case that occurs when all the paths $w_k = x_k w_k x_k$ are
loop edges. Evidently $G^t$ and $G$ are isomorphic, hence the commutant $\flgee^\prime =
\frgee \simeq \fL_{G^t} \simeq \flgee$ is unitarily equivalent to $\flgee$ and is also unitally
partly free. As a variation of this case, instead let $H$ be the subgraph $H = \{ x_k,e_k,w_k :
k\geq 1 \}$. Then $\fL_H$ is unitally partly free, but $\fL_H^\prime \simeq \fL_{H^t}$ is not
even partly free. 


${\mathbf (ii)}$ A non-discrete example is given by the graph $Q$ consisting of vertices $\{ x_q
: q
\in\bbQ\}$ indexed by the rational numbers, and directed edges $e_{qp} = x_q e_{qp} x_p$
whenever $p \leq q$. This example  satisfies the uniform 
aperiodic path property, in fact there is an abundance of infinite non-overlapping directed paths
emanating from each vertex, so
$\fL_Q$ is unitally partly free. However, notice that the quiver algebra $\A_G$ 
is not even partly free (see Section~2.1). Further note that $Q^t$ is 
graph isomorphic to $Q$, thus 
$\fL_Q \simeq \fL_{Q^t} \simeq \fL_Q^\prime$ is unitarily equivalent to its commutant. 


${\mathbf (iii)}$ The following example was suggested to us by Ken Davidson. For $n\geq 2$,
let
$\bbF_n^+$ be the unital free semigroup on $n$ noncommuting letters, written as $\{ 1,2,\ldots,
n\}$, with unit $\phi$. 
Let $G_n$ be the doubly-bifurcating (sideways) infinite tree with vertices $\{ x_w : w
\in\bbF_n^+\}$ indexed by words in $ \bbF_n^+$, and directed edges $\{e_{iw} = x_{iw}
e_{iw} x_w : w \in\bbF_n^+, 1\leq i \leq n \}$ determined by the directions $w\mapsto iw$. 
Then $G_n$ satisfies the uniform aperiodic path property (observe that 
$G_1 = C_\infty$), and
hence $\fL_{G_n}$ is unitally partly free. An interesting point here is that the graph $G_n$
itself has the structure of the full Fock space Hilbert space $\ell^2(\bbF_n^+)$ traced out by its
left creation operators. Whereas, the Fock space $\H_{G_n}$ consists of infinitely
many disjoint infinite-dimensional components, indexed by elements of $\bbF_n^+$. The
transpose graph $G_n^t$ is quite different from $G_n$. In fact the commutant algebra
$\fL_{G_n^t}\simeq \fR_{G_n} = \fL_{G_n}^\prime$ is not partly free. 

${\mathbf (iv)}$ Let $G$ be the directed graph with vertices $\{ x_k : k\geq 1\}$ and edges $\{
e_k = x_k
e_k x_1 : k\geq 1\}$. This graph has  no infinite paths or
double-cycles, hence $\flgee$ is not partly free.
A variation of this example, turning the $e_k$ into non-overlapping paths of length $k$ that only
intersect at vertex $x_1$, produces a non-partly free $\flgee$ with graph containing 
non-overlapping finite paths of arbitrarily large length. 

${\mathbf (v)}$
Let $G= \{ x_1, x_2, e_k=x_2 e_k x_1 : k\geq 1\}$ with $e_k$
distinct edges. Then $\fL_G$ is not partly free and  is unitarily equivalent to its
commutant $\fL_{G}^\prime \simeq \fL_{G^t}$. 


${\mathbf (vi)}$ Let $G$ be the directed graph with vertices $\{ x_k : k\in\bbZ\}$ and directed
edges 
$\{ e_k : k\in\bbZ\}$ where 
\[
e_k = \left\{ \begin{array}{ll}
x_{2m+1} e_{k} x_{2m} & \mbox{if $k=2m$} \\ 
x_{2m+1} e_{k} x_{2m+2} & \mbox{if $k=2m+1$}. 
\end{array}\right. 
\]
Then $\flgee$ is not partly free and is unitarily equivalent to its commutant  $\fL_{G}^\prime
\simeq \fL_{G^t}$.
\end{egs}

\subsection{Partly Free Quiver Algebras}

Using Theorems~\ref{weak} and \ref{strong} we may readily deduce  graph-theoretic conditions
for quiver
algebras $\A_G$ to be partly free. We require the following structural result for partial
isometries in $\A_G$ for the countable graph case. An immediate consequence is that $\A_G$ 
only contains isometries when $G$ has finitely many vertices. 

\begin{lem}\label{quiverlemma}
Let $G$ be a countable directed graph. If $V$ is a partial isometry in $\A_G$, then its initial
projection $V^*V= \sum_{x\in\I} P_x$ is the sum of only finitely many $P_x$.
\end{lem}

\Prf
In fact, if $V$ is a partial isometry in $\flgee$ for which $\I$ is an 
infinite set, then 
$
\dist(V,\A_G) \geq  1. 
$
Indeed, let $q(L)$ belong to the set of polynomials $\A = \Alg \{ L_w : w\in\bbF^+(G) \}$ in the
$L_e$ and $L_x=P_x$. As $\I$ is infinite, there is a $y\in\I$ such that $q(L) P_y = 0$. Hence
$
1 =||VP_y || =  || VP_y - q(L) P_y || \leq || V - q(L) ||,
$
and this proves the claim because $\A$ is (norm) dense in $ \A_G$. The 
lemma follows since $\flgee$ contains $\A_G$. 
\bx

Let $\A_2$ be the quiver algebra generated by the graph with a single vertex and two distinct
loop edges. This is the {\it noncommutative disc algebra} of 
Popescu \cite{Pop_beur,Pop_disc}, also considered by Arias 
\cite{A_disc}, and Muhly-Solel \cite{M1,MS1}.  
Say that $\A_G$ is {\it partly free}, or {\it unitally partly free}, if the maps in
Definition~\ref{partlyfreedefn} are injections of $\A_2$ into $\A_G$, and are restrictions of
injections of the generated $\ca$-algebras. 

We say that a graph $G$ has the {\it double-cycle property} if  $G$ contains 
a 
double-cycle, and  $G$ satisfies the {\it uniform double-cycle 
property}
when
every vertex saturation $\sat(x)$ includes a double-cycle. Compare the 
following results with 
Theorems~\ref{weak} and \ref{strong}, and notice how the proper infinite 
path
phenomena only arises in the $\wot$-closed case. In particular, there are 
many examples in the countable graph case for which $\flgee$ is partly 
free, but $\A_G$ is not. 

\begin{thm}\label{weaknorm}
The following assertions are equivalent for a finite or countable directed graph $G$:
\begin{itemize}
\item[$(i)$] $G$ has the double-cycle property.  
\item[$(ii)$] $\A_G$ is partly free. 
\item[$(iii)$] There are nonzero partial isometries $U$, $V$ in $\A_G$ with
\[
U^*U = V^*V, \,\,\,\, UU^*\leq U^*U, \,\,\,\, VV^*\leq V^*V, \,\,\,\, U^* V =0.
\]
\end{itemize}  
\end{thm}

\Prf
As in Theorem~\ref{weak} it suffices to establish the equivalence of $(i)$ and $(iii)$. But $(i)
\Rightarrow (iii)$ is clear since $\A_G$ will include $U= L_{w_1}$, $V=L_{w_2}$, where
$w_1$, $w_2$ are distinct cycles over a common vertex, when $(i)$ holds. On the
other hand, as $\flgee$ contains $\A_G$, condition $(iii)$ implies $G$ 
satisfies the aperiodic 
path property by Theorem~\ref{weak}. But recall from the proof of $(iii) \Rightarrow (i)$ in
Theorem~\ref{weak} that the proper infinite path part of this property can
only occur when the initial vertex set $\I$ is infinite. Hence by Lemma~\ref{quiverlemma}, $G$
contains a double-cycle and $(i)$ holds. 
\bx

In the unital case, the graph can only have finitely many vertices. 

\begin{thm}\label{strongnorm}
The following assertions are equivalent for a finite or countable directed graph $G$:
\begin{itemize}
\item[$(i)$] $G$ has finitely many vertices and satisfies the 
uniform double-cycle property. 
\item[$(ii)$] $\A_G$ is unitally partly free. 
\item[$(iii)$] There are  isometries $U$, $V$ in $\A_G$ with
\[
U^* V =0.
\]
\end{itemize}  
\end{thm}

\Prf
This result follows from Theorem~\ref{strong} as the previous result follows from
Theorem~\ref{weak}, other
than the extra vertex condition on $G$. In particular, $(iii)$ implies $G$ 
satisfies the uniform 
infinite path property, but since $U$, $V$ are isometries the initial vertex set in this case is the
entire vertex set of $G$, that is $\I = V(G)$. Therefore, $G$ can only have finitely 
many vertices by Lemma~2.11 and
must satisfy the uniform double-cycle property. 
\bx




{\noindent}{\it Acknowledgements.}
We would like to thank Ken Davidson for organizing a workshop on 
nonselfadjoint operator algebras at the Fields Institute in Toronto (July 
2002), where the authors had a number of productive conversations. 
The first named author would also like to thank members of the Department of Mathematics at
Purdue University for kind hospitality during preparation of this article. 

\vskip 1truecm

{\noindent}{\it Note Added in Proof.} 
Problem 2.8 has been answered in the affirmative for all finite graphs \cite{JP}. The problem
remains open for general $\flgee$. We also mention that the ideal structure of $\flgee$ has been
analyzed in \cite{JuryK}. 
Further, Ephrem \cite{Ephrem} has recently identified graph conditions 
which characterize when a Cuntz-Krieger graph $\ca$-algebra has type~I 
representation theory. Interestingly, the condition he 
obtains for finite graphs is equivalent to the graph not having any 
double-cycles.


\baselineskip=12pt

\bigskip

\begin{tabbing}
{\it E-mail address}:xx\= \kill
\noindent {\footnotesize\it Addresses}:
\>{\footnotesize\sc Department of Mathematics and Statistics}\\
\>{\footnotesize\sc University of Guelph}\\
\>{\footnotesize\sc Guelph, Ontario}\\
\>{\footnotesize\sc CANADA \quad N1G 2W1}\\
\\
\>{\footnotesize\sc Department of Mathematics and Statistics}\\
\>{\footnotesize\sc Lancaster University}\\
\>{\footnotesize\sc Lancaster, England}\\
\>{\footnotesize\sc UK \quad LA1 4YW}\\
\\
{\footnotesize\it E-mail addresses}: 
\>{\footnotesize\sf kribs@math.purdue.edu}\\
\>{\footnotesize\sf s.power@lancaster.ac.uk}

\end{tabbing}

\bigskip
          
\bigskip

\thanks{2000 {\it Mathematics Subject Classification.} 47L55, 47L75. }

\end{document}